\documentclass[11pt,twoside]{amsart}
\usepackage{graphicx}
\usepackage{float}
\usepackage{amsmath, amsthm, amscd, amsfonts, amssymb, graphicx, color}
\usepackage[bookmarksnumbered, plainpages]{hyperref}
\newtheorem{thm}{Theorem}[section]

\newtheorem{lem}[thm]{Lemma}

\numberwithin{equation}{section}
\def\pn{\par\noindent}

\begin{document}
\title{GENERAL SUM-CONNECTIVITY INDEX OF TREES AND UNICYCLIC GRAPHS WITH FIXED MAXIMUM DEGREE}
\author{M. K. JAMIL$^*$, I. TOMESCU}

\thanks{{\scriptsize
Received: dd mmmm yyyy, Accepted: dd mmmm yyyy.\\
$*$ Corresponding author}}
\maketitle

\begin{center}
\end{center}

\begin{abstract}
    The general sum-connectivity index of a graph $G$ is defined as $\chi_\alpha(G)=\sum\limits_{uv\in E(G)} {(d(u)+d(v))^{\alpha}}$, where $d(v)$ denotes the degree of the vertex $v$ in $G$ and $\alpha$ is a real number. In this paper it is  deduced the maximum value for the general sum-connectivity index of $n$-vertex trees for $-1.7036\leq \alpha <0$ and of $n$-vertex unicyclic graphs for $-1\le \alpha <0$ respectively, with fixed maximum degree $\triangle $. The corresponding extremal graphs, as well as the $n$-vertex unicyclic graphs with the second maximum general sum-connectivity index for $n\ge 4$ are characterized. This extends the corresponding results by  Du, Zhou and Trinajsti\' c [arXiv.org/1210.5043] about sum-connectivity index.
\end{abstract}

\vskip 0.2 true cm


\pagestyle{myheadings}
\markboth{\rightline {\sl  General sum-connectivity index \hskip 8.5 cm  M. K. JAMIL$^*$, I. TOMESCU}}
         {\leftline{\sl  General sum-connectivity index \hskip 8.5 cm  M. K. JAMIL$^*$, I. TOMESCU}}

\bigskip
\bigskip


\section{Introduction}
Let $G(V(G),E(G))$ be a simple graph, where $V(G)$ and $E(G)$ are the sets of vertices and of edges, respectively. For a vertex $v\in V(G)$, $d(v)$ denotes the degree of vertex $v$, $N(v)$ is the set of vertices adjacent to $v$ and the maximum vertex degree of the  graph $G$ is denoted by $\triangle(G)$. If $uv\in E(G)$, $G-uv$ denotes the subgraph of $G$ obtained by deleting the edge $uv$; similarly is defined the graph $G+uv$ if $uv \notin E(G)$. For $n\ge 3$ let $T(n,\triangle)$ be the set of trees with $n$ vertices and maximum degree $\triangle$ and $U(n,\triangle)$ be the set of unicyclic graphs with $n$ vertices and maximum degree $\triangle$ ($2\le \triangle \le n-1$). Let $P_n$ and $C_n$ be the path and the cycle, respectively, on $n\ge3$ vertices. For $\triangle=2$, $T(n,\triangle)=\{P_n\}$ and $U(n,\triangle)=\{C_n\}$. Attaching a path $P_r$ to a vertex $v$ of a graph means adding an edge between $v$ and a terminal vertex of the path. If $r=1$, then we attach a pendant vertex.\\

The Randi\' c index $R(G)$ was proposed by Randi\' c \cite{r}. This is one of the most used molecular descriptors in structure-property and structure-activity relationship studies \cite{dgo,gf,kh,lg}. The Randi\' c index is also called product-connectivity index and it is defined as
\[R(G)=\sum\limits_{uv\in E(G)} {(d(u)d(v))^{-1/2}}. \]
Bollob\' as and Erd\" os \cite{be} generalized the idea of Randi\' c index and proposed the general Randi\' c index, denoted as $R_\alpha$. It is defined as
\[R_\alpha(G)=\sum\limits_{uv\in E(G)} {(d(u)d(v))^{\alpha}}, \]
where $\alpha$ is a real number.\\

The sum-connectivity index was proposed by Trinajsti\' c et al. \cite{zt} and it was observed that sum-connectivity index and product-connectivity index correlate well among themselves and with the $\Pi$-electronic energy of benzenoid hydrocarbons \cite{ltz}. This concept was extended to the general sum-connectivity index in \cite{zt1} and defined as
\[\chi_{\alpha}(G)=\sum_{uv\in E(G)} {(d(u)+d(v))^{\alpha}}.\]
Then $\chi_{-1/2}$ is the sum-connectivity index \cite{zt}.\\

Several extremal properties of the sum-connectivity and general sum-connectivity indices for trees, unicyclic graphs and general graphs were given in \cite{dz,DZT0,ta,tj,zt,zt1}.\\

In \cite{dzt2} Zhou et al. obtained the maximum sum-connectivity index of graphs in the set of trees and in the set of unicyclic graphs respectively, with a given number of vertices and maximum degree and determined the corresponding extremal graphs. They also found the $n$-vertex unicyclic graphs with the first two maximum sum-connectivity indices for $n\ge4$. In this paper we extend these results for the general sum-connectivity index.\\

\section{Main Results}
First we will discuss two lemmas that will be used in the proofs.

\begin{lem}\cite{DZT0}\label{lem1}
    Let $Q$ be a connected graph with at least two vertices. For $a\ge b\ge 1,$ let $G_1$ be the graph obtained from $Q$ by attaching two paths $P_a$ and $P_b$ to $u\in V(Q)$ and $G_2$ the graph obtained from $Q$ by attaching a path $P_{a+b}$ to $u$. Then $\chi_\alpha(G_2)> \chi_\alpha (G_1)$, for $ \alpha_1\le \alpha <0$, where $\alpha_1\approx -1.7036$ is the unique root of the equation $\frac{3^\alpha-4^\alpha}{4^\alpha-5^\alpha}=2$.
\end{lem}
The following property is an extension of a transformation defined in \cite{dzt2}.
\begin{lem}\cite{dzt2}\label{lem2}
    Let a connected graph $M$ with $|V(M)|\ge 3$ and a vertex $u$ of degree two of $M$. Let $H$ be the graph obtained from $M$ by attaching a path $P_a$ to $u$. Denote by $u_1$ and $u_2$ the two neighbors of $u$ in $M$, and by $u'$ the pendant vertex of the path attached to $u$ in $H$. If $d_H(u_2)\le 3$, then for $H'=H-\{uu_2\}+\{u'u_2\}$ we have $\chi_\alpha(H')> \chi_\alpha (H)$, where $-1\le \alpha< 0$.
\end{lem}
\textbf{Proof.}
     If $d_H(u,u')=1$, then for $\alpha<0$ we have:
    $\chi_\alpha(H')-\chi_\alpha(H)=(d_H(u_1)+2)^\alpha+(d_H(u_2)+2)^\alpha
-(d_H(u_1)+3)^\alpha-(d_H(u_2)+3)^\alpha > 0$.\\
 If $d_H(u,u')\ge 2,$ then
\begin{eqnarray*}
    \chi_\alpha(H')-\chi_\alpha(H)
    &=&(d_H(u_1)+2)^\alpha-(d_H(u_1)+3)^\alpha+(d_H(u_2)+2)^\alpha-(d_H(u_2)+3)^\alpha+2\cdot4^\alpha-3^\alpha-5^\alpha\\
    &>&(d_H(u_2)+2)^\alpha-(d_H(u_2)+3)^\alpha+2\cdot4^\alpha-3^\alpha-5^\alpha.
\end{eqnarray*}
 Since $(x+2)^\alpha-(x+3)^\alpha$ is decreasing for $x\ge 0$, we have $(d_H(u_2)+2)^\alpha-(d_H(u_2)+3)^\alpha\ge 5^\alpha-6^\alpha$. Therefore
$$ \chi_\alpha(H')-\chi_\alpha(H)  > 2\cdot 4^\alpha-3^\alpha-6^\alpha.$$
The function $\eta (x)=2\cdot 4^{x}-3^{x}-6^{x}$ has roots $x_1=-1$ and $x_2= 0$ and $\eta (x)>0$ for $x\in (-1,0)$ \cite{wa}. It follows that $\chi_\alpha(H')>\chi_\alpha(H)$ for every $-1\le \alpha <0$. \hspace{6cm} $\Box$
\\

 For $\frac{n}{2}\le \triangle \le n-1$, let $T_{n,\triangle}$ be the tree obtained by attaching $2\triangle+1-n$ pendant vertices and $n-\triangle-1$ paths of length two to a vertex. For $\frac{n+2}{2} \le \triangle \le n-1$, let $U_{n,\triangle}$ be the unicyclic graph obtained by attaching $2\triangle-n-1$ pendant vertices and $n-\triangle -1$ paths of length two to the same vertex of a triangle.\\

\begin{thm}\label{thm1}
    Let $G\in T(n,\triangle)$, where $2\le \triangle \le n-1$ and $\alpha_1\le \alpha <0$, where $\alpha_1\approx -1.7036$ is the unique root of the equation $\frac{3^\alpha-4^\alpha}{4^\alpha-5^\alpha}=2$. Then
    \[\chi_\alpha(G)\le \left\{ \begin{gathered}
    ((\triangle+2)^\alpha-(\triangle+1)^\alpha+3^\alpha)(n-\triangle-1)+\triangle(\triangle+1)^\alpha \,\, if \,\,\, \,  \frac{n}{2}\le \triangle \le n-1 \hfill \\
    ((\triangle+2)^\alpha+3^\alpha-4^\alpha)\triangle+(n-\triangle-1)4^\alpha\,\,\hspace{2.2cm}  if\,\, 2\le \triangle \le \frac{n-1}{2}\hfill \\
    \end{gathered}  \right.\]
    and equality holds if and only if $G=T_{n,\triangle}$ for $\frac{n}{2}\le \triangle \le n-1$, and $G$ is a tree obtained by attaching $\triangle$ paths of length at least two to a unique vertex for $2\le \triangle \le \frac{n-1}{2}$.
\end{thm}

\textbf{Proof.}
    The case $\triangle=2$ is clear since in this case $G=P_n$. Suppose that $\triangle \ge 3$ and let $G$ be a tree in $T(n,\triangle)$ having maximum general sum-connectivity index. Let $v$ be a vertex of degree $\triangle$ in $G$. If there exists some vertex of degree greater than two in $G$ different from $v$, then by Lemma \ref{lem1}, we may get a tree in $T(n,\triangle)$ with greater general sum-connectivity index, a contradiction. It follows that $v$ is the unique vertex of degree greater than two in $G$. Let $k$ be the number of neighbors of $v$ with degree two. Since in $V(G)\backslash (\{v\}\cup N(v))$
there are $n-\triangle -1$ vertices, it follows that $k\le \min\{n-\triangle-1,\triangle\}$. If $n-\triangle-1\ge \triangle$, i.e., $\triangle\le \frac{n-1}{2}$, then $1\le k\le \triangle$. If $n-\triangle-1<\triangle,$ i.e., $\triangle\ge \frac{n}{2},$ then $0\le k\le n-\triangle-1$.
We get
\begin{eqnarray*}
    \chi_\alpha(G)&=&(\triangle-k)(\triangle+1)^\alpha+k(\triangle+2)^\alpha+k\cdot 3^\alpha+(n-\triangle-k-1)4^\alpha\\
    &=&k((\triangle+2)^\alpha-(\triangle+1)^\alpha+3^\alpha-4^\alpha)+
\triangle(\triangle+1)^\alpha+(n-\triangle-1)4^\alpha.
\end{eqnarray*}
Since $(\triangle+2)^\alpha-(\triangle+1)^\alpha$ is increasing for $\triangle\ge 3$
we obtain $(\triangle+2)^\alpha-(\triangle+1)^\alpha+3^\alpha-4^\alpha\ge 5^\alpha+
3^\alpha-2\cdot4^\alpha>0$, the last inequality holding by Jensen's inequality. Consequently,
    \[\chi_\alpha(G)\le \left\{ \begin{gathered}
    ((\triangle+2)^\alpha-(\triangle+1)^\alpha+3^\alpha)(n-\triangle-1)+\triangle(\triangle+1)^\alpha\,\, if \,\, \frac{n}{2}\le \triangle \le n-1 \hfill \\
    ((\triangle+2)^\alpha+3^\alpha-4^\alpha)\triangle+(n-\triangle-1)4^\alpha\,\, \hskip2.2cm if\,\, 2\le \triangle \le \frac{n-1}{2}.\hfill \\
    \end{gathered}  \right.\]

 For $\frac{n}{2}\le \triangle\le n-1$, the equality holds if and only if $k=n-\triangle-1$, i.e., each of the $n-\triangle-1$ neighbors of degree two  of the vertex $v$  is adjacent to
exactly a pendant vertex, i.e., $G=T_{n,\triangle}$. For  $2\le \triangle\le \frac{n-1}{2}$
the equality holds for $k=\triangle$, i.e., $G$ is a tree obtained by attaching $\triangle$ paths of length at least two to a unique vertex.\hspace{12cm} $\Box$ \\

Now we obtain the maximum general sum-connectivity index of graphs in $U(n,\triangle)$ and deduce the extremal graphs. As a consequence, we deduce the $n$-vertex unicyclic graphs with the first and second maximum general sum-connectivity indices for $n\ge 4$.

\begin{thm}\label{thm2}
    Let $G\in U(n,\triangle),$ where $2\le \triangle \le n-1$ and $-1\le \alpha < 0$. Then
    \[\chi_\alpha(G)\le\left\{\begin{gathered}
    (n-\triangle-1)3^\alpha+(n-\triangle+1)(\triangle+2)^\alpha+(2\triangle-n-1)(\triangle+1)^\alpha+4^\alpha\,\, if\,\, \frac{n+2}{2}\le \triangle \le n-1 \hfill \\
    (\triangle-2)3^\alpha+\triangle(\triangle+2)^\alpha+(n-2\triangle+2)4^\alpha \hskip3cm\,\, if\,\, 2\le\triangle\le \frac{n+1}{2}.
    \end{gathered}\right.\]
    For  $\frac{n+2}{2}\le \triangle \le n-1$ the equality holds if and only if $G=U_{n,\triangle}$. If  $2\le \triangle \le \frac{n+1}{2}$ the equality holds if and only if $G$ is a unicyclic graph obtained by attaching $\triangle-2$ paths of length at least two to a fixed vertex of a cycle.
\end{thm}
\textbf{Proof.}
    The case $\triangle=2$ is trivial since in this case $G=C_n$. Suppose that $\triangle \ge3$, $G$ is a graph in $U(n,\triangle)$ with maximum general sum-connectivity index, and $C$ is the unique cycle of $G$. Let $v$ be a vertex of degree $\triangle$ in $G$.\\
    If $\triangle=3$ and there exists some vertex outside $C$ with degree three, then by Lemma \ref{lem1}, we may get a graph in $U(n,3)$ with greater general sum-connectivity index, a contradiction. If there are at least two vertices on $C$ with degree three, then by Lemma \ref{lem2}, we may deduce the same conclusion. Thus, $v\in V(C)$ and $v$ is the unique vertex in $G$ with degree three. Then either $\chi_\alpha(G)=(n-2)4^\alpha+2\cdot5^\alpha$ when $v$ is adjacent to a vertex of degree one and two vertices of degree two for $n\ge 4$, or $\chi_\alpha(G)=(n-4)4^\alpha+3\cdot5^\alpha+3^\alpha$ when $v$ is adjacent to three vertices of degree two for $n\ge5$. The difference of these two numbers equals
    $(n-2)4^\alpha+2\cdot5^\alpha-(n-4)4^\alpha-3\cdot5^\alpha-3^\alpha=
2\cdot4^\alpha-5^\alpha-3^\alpha<0$ by Jensen's inequality. Hence, $G$ is the graph obtained by attaching a pendant vertex to a triangle for $n=4$, i.e., $G=U_{4,3},$ and a graph obtained by attaching a path of length at least two to a cycle for $n\ge 5$.
\par
 Now suppose that $\triangle\ge4$. As for the case $\triangle =3$ we deduce that the vertex of maximum degree is unique, otherwise $G$ has not a maximum general sum-connectivity index in $U(n,\triangle)$. We will show that the vertex of maximum degree $v$ lies on $C$. Suppose that $v$ is not on $C$. Let $w$ be the vertex on $C$ such that $d_G(v,w)=\min \{d_G(v,x):x\in V(C)\}$. If there is some vertex  outside $C$ with degree greater than two different from $v$, or if there is some vertex on $C$ with degree greater than two different from $w$, then by Lemmas \ref{lem1} and \ref{lem2}, we may get a graph in $U(n,\triangle)$ with greater general sum-connectivity index, a contradiction. Thus, $v$ and $w$ are the only vertices of degree greater than two in $G$, and $d_G(v)=\triangle$ and $d_G(w)=3$. Let $Q$ be the path connecting $v$ and $w$.
Suppose that $v_1,v_2,\cdots,v_{\triangle-1}$ are the neighbors of $v$ outside $Q$. Let $d_i=d_G(v_i)$ for $i=1,\ldots,\triangle-1$. Note that since $G$ has maximum general sum-connectivity index, then $d_1,\ldots ,d_{\triangle -1}\in \{1,2\}$, since otherwise we can apply Lemma \ref{lem1} and obtain a graph having a greater general sum-connectivity index.
 Consider $G_1=G-\{vv_3,\cdots,vv_{\triangle-1}\}+\{wv_3,\cdots,wv_{\triangle-1}\}\in U(n,\triangle )$. Note that $d_{G_1}(w)=\triangle$ and $d_{G_1}(v)=3$. Then
\begin{eqnarray*}
    \chi_\alpha(G_1)-\chi_\alpha(G)&=&(d_1+3)^\alpha-(d_1+\triangle)^\alpha+
(d_2+3)^\alpha-(d_2+\triangle)^\alpha+2(\triangle+2)^\alpha-2\cdot 5^\alpha\\
    &>& 5^\alpha-(2+\triangle)^\alpha+5^\alpha-(2+\triangle)^\alpha+
2(\triangle+2)^\alpha-2\cdot 5^\alpha=0,
\end{eqnarray*}
since the function $(x+3)^{\alpha}-(x+\triangle )^{\alpha}$ is strictly decreasing in $x\geq 0$ for $\triangle \ge 4$.
    Because $d_{G_{1}}(v)=3$, then by Lemma \ref{lem1}, we may get a graph $G'$ in $U(n,\triangle)$ such that $\chi_\alpha(G')>\chi_\alpha(G_1)\ge \chi_\alpha(G)$, a contradiction. Hence, we have shown that $v$ lies on $C$.\\
    If there is some vertex outside $C$ with degree greater than two, then by Lemma \ref{lem1} we may obtain a graph in $U(n,\triangle)$ with greater general sum-connectivity index, a contradiction. If there is some vertex on $C$ with degree three, then by Lemma \ref{lem2}, we may get a graph in $U(n,\triangle)$ with greater general sum-connectivity index, a contradiction. Thus, $G$ is a graph obtained from $C$ by attaching $\triangle-2$ paths to $v$. Let $k$ be the number of neighbors of $v$ with degree two. Then as above we get $k\le \min\{n-\triangle-1,\triangle-2\}$. If $n-\triangle-1\ge \triangle-2$, i.e., $\triangle\le \frac{n+1}{2}$, then $0\le k \le \triangle -2$. If $n-\triangle-1<\triangle-2$, i.e., $\triangle \ge \frac{n+2}{2}$, then $0\le k \le n-\triangle-1$.
 We get
\begin{eqnarray*}
    \chi_\alpha(G)&=&k3^\alpha+(k+2)(\triangle+2)^\alpha+(\triangle-k-2)(\triangle+1)^\alpha+(n-\triangle-k)4^\alpha\\
    &=&(3^\alpha+(\triangle+2)^\alpha-(\triangle+1)^\alpha-4^\alpha)k+(\triangle-2)(\triangle+1)^\alpha+2(\triangle+2)^\alpha+4^{\alpha}(n-\triangle).
\end{eqnarray*}
 We have
$$   3^\alpha-4^\alpha+(\triangle+2)^\alpha-(\triangle+1)^\alpha\ge 3^\alpha-4^\alpha+6^\alpha-5^\alpha >0 $$
since the function $f(x)=(x+2)^{\alpha}-(x+1)^{\alpha}$ is strictly increasing for $\alpha <0$, hence $f(\triangle)\ge f(4)$ and $f(4)>f(2)$. It follows that $\chi _{\alpha}(G)$ is bounded above by
 \[ \left \{\begin{gathered}
    (3^\alpha+(\triangle+2)^\alpha-(\triangle+1)^\alpha-4^\alpha)(n-\triangle-1)+(\triangle-2)(\triangle+1)^\alpha+2(\triangle+2)^\alpha+4^{\alpha}(n-\triangle)\,\\
    \hskip12cmif\,\,\,\frac{n+2}{2}\le \triangle\le n-1 \hfill\\
    (3^\alpha+(\triangle+2)^\alpha-(\triangle+1)^\alpha-4^\alpha)(\triangle-2)+(\triangle-2)(\triangle+1)^\alpha+2(\triangle+2)^\alpha+4^{\alpha}(n-\triangle)\,\, if\,\,\,2\le \triangle\le \frac{n+1}{2} \hfill\\
    \end{gathered}\right.\]

    \[=\left \{\begin{gathered}
    (n-\triangle-1)3^\alpha+(n-\triangle+1)(\triangle+2)^\alpha+(2\triangle-n-1)(\triangle+1)^\alpha+4^\alpha\,\, if\,\,\frac{n+2}{2}\le \triangle\le n-1 \hfill\\
    (\triangle-2)3^\alpha+\triangle(\triangle+2)^\alpha+(n-2\triangle+2)4^\alpha\,\,
\hspace{4.3cm} if \,\,\, 2\le \triangle\le \frac{n+1}{2} \hfill\\
    \end{gathered}\right.\]
Equality holds for $\frac{n+2}{2}\le \triangle \le n-1$ if and only if $k=n-\triangle-1$, i.e., $G=U_{n,\triangle}$; if $2\le \triangle \le \frac{n+1}{2}$ then equality is reached if and only if $k=\triangle-2$, i.e., $G$ is a unicyclic graph obtained by attaching $\triangle-2$ paths of length at least two to a unique vertex of a cycle. \hspace{9cm}$\Box$\\
\begin{thm}
    If $-1\le \alpha <0$, among the unicyclic graphs on $n\ge4$ vertices, $C_n$ is the unique graph with maximum general sum-connectivity index, which is equal to $n4^{\alpha}$. For $n=4$, $U_{4,3}$ is the unique graph with the second maximum general sum-connectivity index, which is equal to  $2\cdot4^\alpha+2\cdot5^\alpha$. For $n\ge5$, the graphs obtained by attaching a path of length at least two to a vertex of a cycle are the unique graphs with the second maximum general sum-connectivity index, which is equal to $(n-4)4^\alpha+3\cdot5^\alpha+3^\alpha$.
\end{thm}

\textbf{Proof.}
    For $n=4$ we get
  \[\chi_\alpha(U_{4,3})-\chi_\alpha(C_4)=2\cdot5^\alpha-2\cdot4^\alpha<0.\]
    Now, suppose that $n\ge5$ and $G$ is a unicyclic graph on $n$ vertices. Let $\triangle$ be the maximum degree of $G$, where $2\le \triangle \le n-1$. Let $f(x)=(x-2)3^\alpha+x(x+2)^\alpha+(n-2x+2)4^\alpha$ for $x\ge 2$. If $\frac{n+2}{2}\le \triangle \le n-1$, then by Theorem \ref{thm2},
    \[\chi_\alpha(G)\le (n-\triangle-1)3^\alpha+(n-\triangle+1)(\triangle+2)^\alpha+(2\triangle-n-1)(\triangle+1)^\alpha+4^\alpha\]
    \[=f(\triangle)+(n-2\triangle+1)(3^\alpha-4^\alpha+(\triangle+2)^\alpha-(\triangle+1)^\alpha)<f(\triangle)\]
since the function $x^{\alpha}-(x+1)^{\alpha}$ is strictly decreasing for $x\ge 0$ and $\triangle \geq 4$.\\
    If $2\le \triangle\le \frac{n+1}{2}$, then by Theorem \ref{thm2}, $\chi_\alpha\le f(\triangle)$ and equality can be reached. We shall prove that $f'(x)<0$, which implies that $f(x)$ is strictly decreasing for $x\geq 2$. One deduces
\[  f'(x)=3^\alpha+(x+2)^\alpha+\alpha x(x+2)^{\alpha-1}-2\cdot4^\alpha. \]
Let
    \[g(x)=(x+2)^\alpha+\alpha x(x+2)^{\alpha-1}.\]
We get
    \[g'(x)=\alpha(x+2)^{\alpha-2}(x(\alpha +1)+4)<0.\]
So, $g(x)$ is strictly decreasing for $x\ge2$, thus implying  $g(x)\le 4^\alpha+2\alpha4^{\alpha-1}$. Consequently,
\[ f'(x)\le 3^\alpha-4^\alpha+2\alpha4^{\alpha-1}=4^{\alpha}\left(\left(\frac{3}{4}\right)^{\alpha}+\frac{\alpha}{2}-1\right).\]
Considering the function $h(x)=\frac{x}{2}+(\frac{3}{4})^x$, we get
    $h''(x)=(ln(\frac{3}{4}))^2(\frac{3}{4})^x>0$,
hence $h(x)$ is strictly convex. Since $h(-1)=5/6<1$ and $h(0)=1$, $h(x)$ being strictly convex on $[-1,0)$, it follows that $h(x)<1$ on this interval, or $f'(x)<0$ for every $-1\le \alpha <0$, hence $f(x)$ is strictly decreasing for $x\ge2$.\\
It follows that for $3<\frac{n+2}{2}\le \triangle\le n-1$ we have $\chi_\alpha(G)
<f(\triangle)<f(3)<f(2)$ and for $3\le \triangle \le \frac{n+1}{2}$ we obtain $\chi_\alpha(G)\le f(\triangle)\le f(3)<f(2)$. It follows that $C_n$ is the unique $n$-vertex unicyclic graph with maximum general sum-connectivity index, equal to $f(2)$. Also the $n$-vertex unicyclic graphs with maximum degree $\triangle =3$ and general sum-connectivity index $f(3)$ are the $n$-vertex graphs with the second maximum general sum-connectivity index. By Theorem \ref{thm2}, these graphs consist from  a cycle $C_l$ of an arbitrary length $l$, $3\le l\le n-2$ and a path of length at least two attached to a vertex of $C_{l}$. \hspace{3cm}$\Box$


%



\bigskip
\bigskip

{\footnotesize \pn{\bf Muhammad Kamran Jamil}\; \\ {Abdus Salam School of Mathematical Sciences},\\ {Government College University,} {Lahore, Pakistan}\\
{\tt Email: m.kamran.sms@gmail.com}\\

{\footnotesize \pn{\bf Ioan Tomescu}\; \\ {Faculty of Mathematics and Computer Science, {University
of Bucharest,} \\{Str. Academiei, 14, 010014} {Bucharest, Romania.}\\
{\tt Email: ioan@fmi.unibuc.ro}\\

\end{document}